\documentclass{article}
\usepackage{amssymb}
\usepackage{amsmath}
\usepackage{amsfonts}

\newtheorem{theorem}{Theorem}

\newtheorem{definition}[theorem]{Definition}

\newtheorem{lemma}[theorem]{Lemma}

\newtheorem{proposition}[theorem]{Proposition}

\input{tcilatex}
\begin{document}

\title{The canonical geometry\\
of a Lie group}
\author{Erc\"{u}ment H. Orta\c{c}gil}
\maketitle

\begin{abstract}
An abstract Lie group $G$ admits many left-invariant metrics and it is well
known that these metrics posess drastically different curvature properties.
However, $G$ admists a \textit{canonical }metric if we view $G$ as a flat
and globalizable absolute parallelism $w$ as in [O1]. We study some
surprising consequences of this shift in perspective.
\end{abstract}

\section{Basic concepts}

Our main object is a pair $(M,w)$ where $M$ is a smooth manifold and $%
w=(w_{j}^{i}(x))$ is a geometric object on $M,$ called the structure object
in [O1] and [O2]. To understand the meaning of $w,$ let $F(M)\rightarrow M$
be the principal frame bundle of $M$ whose fiber over $p\in M$ is the set of 
$1$-jets (which we call $1$-arrows) $j_{1}(f)^{o}$ of local diffeomorphisms $%
f$ with the source at the origin $o\in \mathbb{R}^{n}$ and the target at $%
f(o)=p.$ The $1$-arrows with the source and target at $o$ is the Lie group $%
GL(n,\mathbb{R}),$ $n=\dim M,$ and acts freely on the fiber over $p$ by
composition at the source. Therefore $F(M)\rightarrow M$ is a right
principal bundle with the structure group $GL(n,\mathbb{R}).$ Now suppose
that $F(M)\rightarrow M$ is trivial and we fix one trivialization $w$ once
and for all. For a coordinate neighborhood $(U,x)\subset M,$ the unique $1$%
-arrow $j_{1}(f)^{o}$ with target at $f(o)=x\in (U,x)$ is of the form $%
w_{j}^{i}(x)=\left[ \frac{\partial f^{i}(z)}{\partial z^{j}}\right] _{z=o}$
where $(z^{i})$ are the standard coordinates in $\mathbb{R}^{n}.$ By the
chain rule, a coordinate change $(U,x)\rightarrow (U,y)$ at the target
transforms the components $(w_{j}^{i}(x))$ acording to

\begin{equation}
w_{j}^{i}(y)=\frac{\partial y^{i}}{\partial x^{a}}w_{j}^{a}(x)
\end{equation}%
Therefore, a trivialization determines a geometric $w$ on $M$ with
components $(w_{j}^{i}(x))$ on $(U,x)$ subject to (1). Conversely, the
geometric object $w=(w_{j}^{i}(x))$ subject to (1) defines uniquely a
trivialization in the obvious way. Therefore, trivializations and the
geometric objects $w$ on $M$ satisfying (1) can be identified and henceforth
we will adhere to this identification but remark here that the geometric
object viewpoint generalizes in a natural way to geometric structures other
than trivializations. Note that if $\widetilde{w}=w^{-1},$ i.e., $\widetilde{%
w}_{a}^{i}(x)w_{j}^{a}(x)=w_{a}^{i}(x)\widetilde{w}_{j}^{a}(x)=\delta
_{j}^{i},$ then $\widetilde{w}=(\widetilde{w}_{j}^{i}(x))$ is another
trivialization subject to

\begin{equation}
\widetilde{w}_{j}^{i}(y)=\widetilde{w}_{a}^{i}(x)\frac{\partial x^{a}}{%
\partial y^{j}}
\end{equation}%
and (2) is obtained by inverting (1). Clearly, $\widetilde{w}$ can be
interpreted as the corresponding trivialization of the coframe bundle $%
\widetilde{F}(M)\rightarrow M$ whose fiber over $p$ is the set of $1$-arrows
with source at $p$ and target at $o\in \mathbb{R}^{n}$ and $\widetilde{F}%
(M)\rightarrow M$ is a left $GL(n,\mathbb{R})$-principal bundle.

Given $(M,w),$ there are three objects canonically associated with $w$ which
will play a fundamental role below. The first is a very special groupoid $%
\Upsilon $ on $M.$ For $p,q\in M,$ we define a $1$-arrow $\varepsilon (p,q)$
from $p$ to $q$ by composing the inverse of the $1$-arrow from $o$ to $p$
with the $1$-arrow from $o$ to $q.$ This amounts to defining $\varepsilon
(p,q)$ in coordinates as

\begin{equation}
\varepsilon _{j}^{i}(x,y)=w_{a}^{i}(y)\widetilde{w}_{j}^{a}(x)
\end{equation}

We easily check that the $1$-arrows defined by $\varepsilon $ are closed
under composition and inversion of arrows and we obtain a groupoid $\Upsilon 
$ as a subgroupoid of the universal jet groupoid $\mathcal{U}_{1}.$ By
construction, $w$ is left invariant by the arrows of $\Upsilon ,$ i.e., $%
\Upsilon \subset \mathcal{U}_{1}$ is the invariance subgroupoid of $w.$ Part
1 of [O1] and the second chapter of [O2] are devoted to a detailed study of $%
\Upsilon $ which we will assume henceforth. It is crucial to observe that
different $w$ may define the same groupoid $\Upsilon .$ Indeed, let $A\in
GL(n,\mathbb{R)}$ and $(M,w)$ be given. We define $Aw$ by $\left( Aw\right)
_{j}^{i}(x)\overset{def}{=}A_{j}^{a}w_{a}^{i}(x),$ that is, $Aw$ is the
trivialization obtained by acting on the $1$-arrows of $w$ at the\textit{\
source }with the constant matrix $A.$ However, note that $A$ cancels in (3)
and $w$ and $Aw$ define the same groupoid $G.$

The second object naturally associated with $w$ is the canonical metric $g$
on $M$ whose value $g(p)$ at $p\in M$ is obtained by mapping the standard
Euclidean metric of $\mathbb{R}^{n}$ to $p\in M$ using the $1$-arrow of the
trivialization $w$ from $o$ to $p.$ In coordinates, we have

\begin{equation}
g^{ij}(x)\overset{def}{=}\sum_{1\leq a\leq n}w_{a}^{i}(x)w_{a}^{j}(x)\text{
\ \ \ \ }g_{ij}(x)=\sum_{1\leq a\leq n}\widetilde{w}_{i}^{a}(x)\widetilde{w}%
_{j}^{a}(x)\text{ \ \ \ \ \ \ \ }\widetilde{w}=w^{-1}
\end{equation}

Clearly $g=(g_{ij}(x))$ is symmetric, positive definite and also $\Upsilon $%
-invariant. Now if we fix the trivialization $w,$ then it is easily checked
that the metrics of $Aw$ range over all $\Upsilon $-invariant metrics on $M$
as $A$ ranges over $GL(n,\mathbb{R}).$ If $\mathfrak{R}=0$ and the groupoid $%
\Upsilon $ integrates to a globalizable $\mathcal{G},$ i.e., if the
pseudogroup $\mathcal{G}$ of local solutions globalizes to a transitive
transformation group $\mathcal{G}$ of $M$ acting simply transitively on $M,$
then these metrics are of course also $\mathcal{G}$-invariant. Using $%
\mathcal{G},$ we can now define a Lie group structure on $M$ in such a way
that $\mathcal{G}$ becomes left or right translations acording to our choice
and therefore all these metrics become left or right invariant. As a
surprising fact, it turns out that these metrics have drastically different
curvature properties (see the excellent survey article [M] and the recent
book [AB]). However, if we shift our focus from metrics to transformations,
i.e., from Riemannian geometry to Lie theory as proposed in [O1] where an
abstract Lie group is by definition a flat and globalizable $w$ modulo the
choices of left/right and a unit, then the primary object for us is $w$
which has the canonical metric (4) ! Though it does not interest here, a
trivialization $(M,w)$ has also a canonical syplectic form, almost complex
structure...etc...any first order structure which has a canonical meaning in 
$\mathbb{R}^{n}$ carries over $M$ by $w.$

Finally, the third object naturally associated with $w$ is the integrability
object $I=(I_{jk}^{i}(x))$ defined by

\begin{equation}
I_{jk}^{i}(x)\overset{def}{=}\left[ \frac{\partial w_{a}^{i}(x)}{\partial
x^{j}}\widetilde{w}_{k}^{a}(x)\right] _{[jk]}=\frac{\partial w_{a}^{i}(x)}{%
\partial x^{j}}\widetilde{w}_{k}^{a}(x)-\frac{\partial w_{a}^{i}(x)}{%
\partial x^{k}}\widetilde{w}_{j}^{a}(x)
\end{equation}

Both the linear and nonlinear curvatures $\mathfrak{R}$ and $\mathcal{R}$ in
[O1] are determined by $I.$ Now $w$ is $\Upsilon $-invariant by definition
but $I$ is not necessarily $\Upsilon $-invariant. It turns out that $%
\mathfrak{R}=0\iff \mathcal{R}=0\iff I$ is $\Upsilon $-invariant ([O1],
[O2]).

In Part 1 of [O1] we introduced two linear connections $\widetilde{\nabla },$
$\nabla $ both determined by $w$ but in [O2] we showed that these
connections are consequences of more fundamental concepts: $\widetilde{%
\nabla }$ is formal Lie derivative and $\nabla $ is the Spencer operator.
Therefore, it is possible to develop the whole theory without even
mentioning the word "linear connection" and this approach applies to all
geometric structures ([O2]). However, since geometers are much familiar with
linear connections, we will continue to use this interpretation below.

\section{The curvature of a local Lie group}

We start with the first Bianchi identity (Proposition 6.4 in [O1])

\begin{equation}
\mathfrak{R}_{kj,r}^{i}+\mathfrak{R}_{jr,k}^{i}+\mathfrak{R}%
_{rk,j}^{i}=I_{kj}^{a}I_{ar}^{i}+I_{jr}^{a}I_{ak}^{i}+I_{rk}^{a}I_{aj}^{i}
\end{equation}

We recall that $\mathfrak{R}$ is the linear curvature of $(M,w)$ defined by $%
\widetilde{\nabla }_{r}I_{kj}^{i}=\mathfrak{R}_{kj,r}^{i}$ in [O1] (a more
conceptual definition is given in [O2])) and it turns out that $I$ and $%
\mathfrak{R}$ are the torsion and curvature of the linear connection $\nabla
.$ We define

\begin{equation}
\mathcal{S}_{kj,r}^{i}\overset{def}{=}I_{kj}^{a}I_{ar}^{i}
\end{equation}

We observe that we differentiate $w$ twice to define $\mathfrak{R}$ whereas
we differentiate $w$ only once to define $\mathcal{S}.$ This suggests that $%
\mathcal{S}$ is a more fundamental object than $\mathfrak{R}.$

\begin{definition}
$\mathcal{S}$ is the primary curvature of $(M,w).$
\end{definition}

Obviously $\mathcal{S}$ satisfies%
\begin{equation}
\mathcal{S}_{kj,r}^{i}=-\mathcal{S}_{jk,r}^{i}
\end{equation}

Now suppose $(M,w)$ is a local Lie group (LLG), i.e., $\mathfrak{R}=0.$ In
this case we denote $(M,w)$ by $(M,w,\mathcal{G)}$ where $\mathcal{G}$ is
the pseudogroup obtained by locally integrating the $1$-arrows of the
groupoid $\Upsilon .$ Now the left hand side (LHS) of (6) vanishes and $%
\mathcal{S}$ satisfies also

\begin{equation}
\mathcal{S}_{kj,r}^{i}+\mathcal{S}_{jr,k}^{i}+\mathcal{S}_{rk,j}^{i}=0
\end{equation}

We define

\begin{equation}
\mathcal{S}_{kj,ri}\overset{def}{=}\mathcal{S}_{kj,r}^{a}g_{ai}
\end{equation}

The main idea of this note is based on the following elementary observation.

\begin{lemma}
On a LLG $(M,w,\mathcal{G)},$ we also have

\begin{equation}
\mathcal{S}_{kj,ri}=-\mathcal{S}_{kj,ir}
\end{equation}
\end{lemma}

\bigskip Since the proof we will make use $\nabla ,\widetilde{\nabla },$ it
is useful to use a notation coherent with these operators. So we make the
definitions%
\begin{equation}
\Gamma _{jk}^{i}\overset{def}{=}\frac{\partial w_{a}^{i}(x)}{\partial x^{j}}%
\widetilde{w}_{k}^{a}(x)\text{ \ \ \ \ \ \ \ \ \ \ }T_{jk}^{i}\overset{def}{=%
}\Gamma _{jk}^{i}-\Gamma _{kj}^{i}=I_{jk}^{i}
\end{equation}

Note that $\Gamma =(\Gamma _{jk}^{i})$ are the components of $\nabla $ and
not the Christoffel symbols of the canonical metric $g.$ In particular $%
\Gamma _{jk}^{i}$ is not necessarily symmetric in $j,k.$ As we remarked
above, $I=T=(T_{jk}^{i})$ and $\mathfrak{R=(R}_{jk,r}^{i}\mathfrak{)}$ are
the torsion and curvature of $\nabla .$ Now a straightforward computation
using the definitions gives the formula (called the Ricci formula in the old
books)

\begin{eqnarray}
\nabla _{s}\nabla _{r}g_{ij}-\nabla _{r}\nabla _{s}g_{ij} &=&-\mathfrak{R}%
_{sr,ij}-\mathfrak{R}_{sr,ji}-T_{sr}^{a}\nabla _{a}g_{ij} \\
&=&-T_{sr}^{a}\nabla _{a}g_{ij}  \notag
\end{eqnarray}

We also have

\begin{eqnarray}
\nabla _{r}g_{ij} &=&\widetilde{\nabla }%
_{r}g_{ij}+T_{ir}^{a}g_{aj}+T_{jr}^{a}g_{ai} \\
&=&T_{ir}^{a}g_{aj}+T_{jr}^{a}g_{ai}  \notag
\end{eqnarray}%
since $g$ is $\Upsilon $-invariant and therefore $\widetilde{\nabla }$%
-parallel. Substituting (14) into (13) gives

\begin{equation}
\nabla _{s}\nabla _{r}g_{ij}-\nabla _{r}\nabla
_{s}g_{ij}=-T_{sr}^{a}T_{ia}^{b}g_{bj}-T_{sr}^{a}T_{ja}^{b}g_{bi}
\end{equation}

Our purpose is to show that the RHS of (15) vanishes. Now we will compute
the LHS of (15) in a different way. Applying $\nabla _{s}$ to (14) gives

\begin{equation}
\nabla _{s}\nabla _{r}g_{ij}=(\nabla _{s}T_{ir}^{a})g_{aj}+T_{ir}^{a}\nabla
_{s}g_{aj}+(\nabla _{s}T_{jr}^{a})g_{ai}+T_{jr}^{a}\nabla _{s}g_{ai}
\end{equation}

We have

\begin{equation}
\nabla _{s}T_{ir}^{a}=\widetilde{\nabla }%
_{s}T_{ir}^{a}+T_{sa}T_{ir}^{a}+T_{ra}T_{si}^{a}+T_{ia}T_{rs}^{a}
\end{equation}

Now we substitute (17) and (14) into (16), alternate $s,r$ in (16) and
equate it to (15). Simplifying the resulting equality using the Bianchi
identity (6) and $\widetilde{\nabla }_{r}T_{jk}^{i}=\widetilde{\nabla }%
_{r}I_{jk}^{i}=\mathfrak{R}_{jk,r}^{i}=0,$ after a straightforward
computation we obtain the desired result

\begin{equation}
T_{sr}^{a}T_{aj}^{b}g_{bi}=-T_{sr}^{a}T_{ai}^{b}g_{bj}
\end{equation}

To summarize what we have done so far, we state

\begin{proposition}
On a LLG $(M,w,\mathcal{G}),$ the primary curvature $\mathcal{S}$ satisfies
the identities
\end{proposition}

\begin{eqnarray}
\mathcal{S}_{ij,kr} &=&-\mathcal{S}_{ij,kr} \\
\mathcal{S}_{ij,kr}+\mathcal{S}_{jk,ir}+\mathcal{S}_{ki,jr} &=&0  \notag \\
\mathcal{S}_{ij,kr} &=&-\mathcal{S}_{ij,rk}  \notag
\end{eqnarray}

Therefore $\mathcal{S}$ \textit{must }satisfy also

\begin{equation}
\mathcal{S}_{ij,kr}=\mathcal{S}_{kr,ij}
\end{equation}%
(see [KN] for a coordinate free derivation of this fact).

From (10) and (11) we deduce

\begin{equation}
\mathcal{S}_{ij,k}^{r}=-\mathcal{S}_{ij,r}^{k}
\end{equation}%
and therefore

\begin{equation}
\mathcal{S}_{ij,a}^{a}=0
\end{equation}

Now summing over $i,r$ in (6) and using (22), we get 
\begin{equation}
\mathcal{S}_{ak,j}^{a}=\mathcal{S}_{aj,k}^{a}
\end{equation}

We define

\begin{equation}
Ric(\mathcal{S})_{kj}\overset{def}{=}\mathcal{S}_{ak,j}^{a}
\end{equation}

\begin{equation}
Ric(\mathcal{S})_{a}^{a}=Ric(\mathcal{S})_{ba}g^{ab}\overset{def}{=}K
\end{equation}%
and note that $Ric(\mathcal{S})_{kj}$ is symmetric by (23).

\begin{definition}
$Ric(\mathcal{S})$ is the Ricci curvature and $K$ is the scalar curvature of
the LLG $(M,w,\mathcal{G}).$
\end{definition}

At this point, it is natural to suspect that $\mathcal{S}$ is the Riemann
curvature tensor $R$ of the canonical metric $g.$ Clearly, without the
assumption $\mathfrak{R}=0,$ there is no reason to believe this and indeed
this belief is not justified. From our stanpoint of Lie theory, however, it
is more natural to ask first the geometric meaning of $\mathcal{S}$ \textit{%
in terms of the group structure of the LLG }$(M,w,\mathcal{G}).$ For this
purpose, we first take a closer look at (14). In addition to $\widetilde{%
\nabla }g=0,$ suppose we also have $\nabla g=0.$ This gives

\begin{equation}
T_{ri}^{a}g_{aj}=-T_{rj}^{a}g_{ai}
\end{equation}

Multiplying (26) with $T_{sm}^{r}$ and summing over $r$ gives

\begin{equation}
T_{sm}^{b}T_{bi}^{a}g_{aj}=-T_{sm}^{b}T_{bj}^{a}g_{ai}
\end{equation}%
which is (11). Therefore (26) is a stronger condition than (11). If the
pseudogroup $\mathcal{G}$ globalizes, then the conditions $\widetilde{\nabla 
}g=\nabla g=0$ are equivalent to the bi-invariance of the metric $g$ with
respect to $\mathcal{G}$ and its centralizer $C(\mathcal{G}).$ Thus we see
that Lemma 2 drops the hypotheses of globalizability and bi-invariance, but
makes the weaker conclusion (27) rather than (26).

Now we inspect the first formula of (6) more closely. The index $a$ in $%
w_{a}^{i}(x)w_{a}^{j}(x)$ represents $\left[ \frac{\partial }{\partial y^{a}}%
\right] _{y=o}$ where $y=(y^{i})$ are the standard coordinates in $\mathbb{R}%
^{n}.$ According to (1), $w_{a}^{i}(x)$ does not transform in the index $a$
and tranforms like a vector field in the index $i.$ Identifying the
variables $y^{1},...,y^{n}$ with $1,2,...,n,$ we now define $n$ \textit{%
global }vector fields $w_{(k)}$ on $M$ by

\begin{equation}
w_{(k)}=(w_{(k)}^{i}(x))=w_{(k)}^{a}\frac{\partial }{\partial x^{a}}\text{ \
\ \ \ \ \ }1\leq k\leq n
\end{equation}

According to (4) and (2), we have

\begin{eqnarray}
g_{ab}w_{(k)}^{a}w_{(l)}^{b} &=&\left( \sum_{1\leq c\leq n}\widetilde{w}%
_{a}^{(c)}(x)\widetilde{w}_{b}^{(c)}(x)\text{ }\right) w_{(k)}^{a}w_{(l)}^{b}
\\
&=&\sum_{1\leq c\leq n}\delta _{(k)}^{(c)}\delta _{(l)}^{(c)}=\delta
_{(k),(l)}  \notag
\end{eqnarray}

Therefore the vector fields $w_{(k)}$ are orthonormal. Since $\widetilde{%
\nabla }g=0,$ they are also $\Upsilon $-invariant and form $n$ orthonormal
global vector fields on $M.$ Note that the assumption $\mathfrak{R}=0$
localizes any object $A$ with $\widetilde{\nabla }A=0$ at some arbitrary
point and therefore reduces all computations to pure algebra. Clearly, if $%
\mathfrak{R}=0,$ then $\widetilde{\nabla }\mathcal{S}=0$ and therefore also $%
\widetilde{\nabla }Ric(\mathcal{S})=0$ and $\widetilde{\nabla }K=0.$
Henceforth we will denote $w_{(k)}$ by $w_{k}.$

Now we define

\begin{equation}
S_{kl}\overset{def}{=}-\mathcal{S}%
_{ab,cd}w_{k}^{a}w_{l}^{b}w_{k}^{c}w_{l}^{d}
\end{equation}

The minus sign will be clear shortly. The RHS of (30) has a coordinatefree
meaning since we contract tensors. We observe that $S_{kl}$ are \textit{not
the components of a tensor} but are numbers depending on $w_{k},$ $w_{l}$
which are defined \textit{canonically} on $M.$ Since $\widetilde{\nabla }%
w_{k}=0,$ $1\leq k\leq n,$ and also $\widetilde{\nabla }\mathcal{S}=0$ if $%
\mathfrak{R}=0,$ on a LLG $(M,w,\mathcal{G})$ we have $S_{kl}(p)=S_{kl}(q)$
for all $p,q\in M,$ i.e., $S_{kl}$ is constant on $M$ for all $1\leq k,l\leq
n.$ However, note that $S_{kl}$ and $S_{jm}$ need not be the same constants.
Clearly $K$ is also constant on $(M,w,\mathcal{G)}.$

\begin{definition}
$S_{kl}$ is the sectional curvature determined by $w_{k}$ and $w_{l}.$
\end{definition}

For the moment, $S$ is a function which assigns numbers to the pairs $w_{k},$
$w_{l}$ but it is easy to extend its definition to \textit{all planes}. A
miracle is hidden in (30): Summing over $l$ and using (4), we get

\begin{eqnarray}
\sum_{1\leq l\leq n}S_{kl} &=&-\mathcal{S}_{ab,cd}w_{k}^{a}w_{k}^{c}g^{bd} \\
&=&-\mathcal{S}_{ab,c}^{b}w_{k}^{a}w_{k}^{c}  \notag \\
&=&\mathcal{S}_{ba,c}^{b}w_{k}^{a}w_{k}^{c}  \notag \\
&=&Ric(\mathcal{S})_{ab}w_{k}^{a}w_{k}^{b}  \notag
\end{eqnarray}

Summing over $k$ in (31) gives

\begin{equation}
\sum_{1\leq l,k\leq n}S_{kl}=\sum_{1\leq k\leq n}Ric(\mathcal{S}%
)_{ab}w_{k}^{a}w_{k}^{b}=Ric(\mathcal{S})_{ab}g^{ab}=Ric(\mathcal{S}%
)_{a}^{a}=K
\end{equation}

Turning back to the question of the group theoretic meaning of $\mathcal{S},$
we fix a point $p\in M$ and consider all paths $c(t)$ with $c(0)=p$
satisfying the following condition: The translation of the tangent vector $%
\overset{\bullet }{c}(0)$ to $c(t)$ by the $1$-arrow of $\Upsilon $ from $%
c(0)$ to $c(t)$ gives the tangent vector $\overset{\bullet }{c}(t)$ (see
[O1], pg.51). These paths are defined for all $t$ and we call them pre-$1$%
-parameter subgroups. If $\mathfrak{R}=0$ and $\mathcal{G}$ globalizes to a
Lie group, they become $1$-parameter subgroups. Differentiation of the pre-$%
1 $-parameter condition gives the geodesics of the linear connections $%
\widetilde{\nabla }$ and $\nabla .$ Therefore, from our standpoint, pre-$1$%
-parameter condition is more fundamental than the geodesic condition. Now we
fix a $2$-dimensional subspace of $T_{p}M.$ As the tangent vectors range
over this $2$-dimensional subspace, the pre-$1$-parameter subgroups
spreading out from $p$ with those tangents sweep a surface. Now suppose $%
\mathfrak{R}=0$ and assume for simplicity that $\mathcal{G}$ globalizes to a
Lie group. Do the transformations of $\mathcal{G}$ permute these surfaces?
Equivalently, for $p,q\in M,$ does the unique transformation of $\mathcal{G}$
that maps $p$ to $q$ restrict to the slices passing through $p,q?$ If so, we
get a $2$-dimensional foliation on $M$ whose slices are subgroups of $%
\mathcal{G}.$ However, even though a Lie group always has $1$-dimensional
subgroups, it does not always admit $2$-dimensional subgroups, i.,e, it is
not always possible to put together $1$-parameter subgroups and form higher
dimensional subgroups. However, this can be done if $\mathcal{G}$ is
solvable. Recalling the Cartan-Killing form and its interpretation in terms
of solvability, it is natural to suspect a relation between vanishing of $%
\mathcal{S}$ and solvability. So we now take a closer at (7). We recall that
the algebraic bracket of two vector fields $\xi ,\eta $ on $M$ is defined by

\begin{equation}
\{\xi ,\eta \}=T_{ab}^{i}\xi ^{a}\eta ^{b}=I_{ab}^{i}\xi ^{a}\eta ^{b}
\end{equation}

Note that $\{,\}$ is an algebraic operation and does not satisfy the Jacobi
identity in general. However, if $\xi ,\eta $ are $\Upsilon $-invariant and $%
\mathfrak{R}=0,$ then $\{\xi ,\eta \}$ coincides with their true bracket $%
[\xi ,\eta ].$ From (7) and (33), we deduce

\begin{equation}
\mathcal{S(}\xi ,\eta ;\gamma )=\mathcal{S}_{bc,d}^{i}\xi ^{b}\eta
^{c}\gamma ^{d}=I_{bc}^{a}I_{ar}^{i}\xi ^{b}\eta ^{c}\eta ^{d}=\{\{\xi ,\eta
\},\gamma \}
\end{equation}

Now we assume $\mathfrak{R}=0$ so that

\begin{equation}
\mathcal{S(}\xi ,\eta ;\gamma )=[[\xi ,\eta ],\gamma ]
\end{equation}%
where $\xi ,\eta ,\gamma $ are $\Upsilon $-invariant, i.e., they belong to
the "Lie algebra" of $\mathcal{G}.$ Denoting the Killing form by $\varkappa
, $ we have

\begin{eqnarray}
\varkappa (\xi ,\eta ) &=&Tr(ad(\xi )\circ ad(\eta ))=Tr\left( x\rightarrow
\lbrack \xi ,[\eta ,x]]\right) \\
&=&Tr(x\rightarrow \lbrack \lbrack x,\eta ],\xi ])=Tr(x\rightarrow \mathcal{%
S(}x,\eta ;\xi ))  \notag \\
&=&Tr(x^{i}\rightarrow \mathcal{S}_{ab,c}^{i}x^{a}\xi ^{b}\eta ^{c})=%
\mathcal{S}_{ab,c}^{a}\xi ^{b}\eta ^{c}  \notag \\
&=&Ric(\mathcal{S)(}\xi ,\eta )  \notag
\end{eqnarray}

We call the LLG $(M,w,\mathcal{G})$ semisimple (solvable) if the Lie algebra 
$\mathfrak{G}$ of the $\Upsilon $-invariant vector fields is semisimple
(solvable). From (35), we see that $\mathcal{S}=0$ if and only if $\mathfrak{%
G}^{(2)}=[[\mathfrak{G,G],G]=}0,$ i.e., $\mathfrak{G}$ is $2$-step
nilpotent. Combining this fact with the Cartan-Killing criterion, we obtain

\begin{proposition}
Let $(M,w,\mathcal{G})$ be a LLG. Then

1) $\mathcal{S}=0$ if and only if $(M,w,\mathcal{G})$ is $2$-step nilpotent

2) If $(M,w,\mathcal{G})$ is nilpotent, then $Ric(\mathcal{S})=0$

3) $(M,w,\mathcal{G})$ is semisimple if and only if $Ric(\mathcal{S})$ is
nondegeretate
\end{proposition}

It is known that a compact Lie group admits a bi-invariant metric which is
not necessarily unique. The existence and uniqueness of a bi-invariant
metric on an abstract Lie group is a subtle issue (see [M], [AB]). On the
other hand, we observe that Proposition 6 assumes neither compactness, nor
globalizability nor the existence of any extra structure like a bi-invariant
metric and it is worthwhile to compare its conceptual simplicity to [M],
[AB]. Some other interesting results follow from the present framework. For
instance, using the linear Spencer sequence with representations we can
define closed forms in the Lie algebra cohomology using $\mathcal{S},$ we
can refine, extend and reinterpret the results in [YB], [G] about the Betti
numbers in terms of the group theory of $(M,w,\mathcal{G})...$We will
relegate these problems to some mythical future work.

Now the following proposition should not come as a surprise.

\begin{proposition}
$-\mathcal{S=}R=$ the Riemann curvature tensor of the canonical metric $g.$
\end{proposition}

The proof is not difficult and we will leave it to the interested reader
(The minus sign suggests to modify (7) by sign, see (40) below). Therefore,
by a classical result (see [C]), the above sectional curvatures are the
Gaussian curvatures of the sweeping surfaces, indicating remarkable
relations between metric and Lie theoretic properties of $(M,w,\mathcal{G}).$

Now given the trivialization $(M,w),$ we define

\begin{equation}
R_{jk}^{i},_{r}\overset{def}{=}\mathfrak{R}_{kj,r}^{i}-I_{kj}^{a}I_{ar}^{i}
\end{equation}

Clearly we have

\begin{equation}
R_{jk}^{i},_{r}=-R_{kj}^{i},_{r}
\end{equation}

\begin{equation}
R_{jk}^{i},_{r}+R_{kr}^{i},_{j}+R_{rj}^{i},_{k}=0
\end{equation}

We will write (37) in the form

\begin{equation}
R=\mathfrak{R-}\mathcal{S}
\end{equation}

To close the scene, we will state here the following \textit{Decomposition
Theorem for absolute parallelism }which we hope will attract the attention
of young researchers and add more suspense to this adventure .

\begin{theorem}
$R$ is the Riemann curvature tensor $R$ of the canonical metric of the
trivialization $(M,w)$ and decomposes as (40).
\end{theorem}

Note that (40) decomposes a metrical object into two pre-Lie theoretic
objects. If $\mathfrak{R}=0,$ then $R$ becomes $-\mathcal{S}.$ If $\mathcal{S%
}$ vanishes, then the RHS of (6) vanishes, i.e., $\{,\}$ satisfies the
Jacobi identity. This condition does not imply $\mathfrak{R}=0$ but makes
the trivialization $(M,w)$ in some vague sense "close to a $2$-step
nilpotent LLG", i.e., almost $2$-step nilpotent by Proposition 6. With this
assumption, we observe the intriguing fact that the homogeneous tensor $%
\mathcal{H=(H}_{j}^{i}\mathcal{)}$ defined in Chapter 13 of [O1] becomes

\begin{equation}
\mathcal{H}_{j}^{i}=-Ric(R)_{(j)}^{i}
\end{equation}%
with the passive index $(j)$ as explained in [O1].

Theorem 8 has many consequences and here is an immediate one which is
already nontrivial: Can two Lie groups $G_{1},$ $G_{2}$ with their left
invariant metrics $g_{1},g_{2}$ be globally isometric but nonisomorphic as
Lie groups? The answer is affirmative and such examples abound everywhere.
For instance, let $G$ be a simply connected $2$-step nilpotent Lie group
which is nonabelian. By (40) it is globally isometric to $\mathbb{R}^{n},$ $%
\dim G=n,$ but $\mathbb{R}^{n}$ is also abelian with its left invariant
Euclidean metric. This example shows that a "Lie structure" is much more
refined than a "metric structure". Understanding this relation, we believe,
is equivalent to understanding the meaning of (40), i.e., the meaning of the
first Bianchi Identity (6) (see Proposition 12.2, [O1] for the second BI).

A final remark: A trivial principle bundle is a very boring object for a
topologist. It is a great wonder that a trivialization of the principal
frame bundle (a first order jet bundle) is such an immensely rich geometric
structure.

\bigskip

\bigskip

\textbf{References}

\bigskip

[AB] Alexandrino, M.M., Bettiol, R.G., Lie Groups and Geometric Aspects of
Isometric Actions, Springer, 2015

[C] do Carmo, M.P., Riemannian Geometry, Birkhauser, 1992

[G] Goldberg, S.I., Curvature and Homology, Academic Press, New York, 1962

[KN] Kobayashi, S., Nomizu, K., Foundations of Differential Geometry,
Interscience Publishers, John Wiley \& Sons, Vol1, 1963

[M] Milnor, J,W., Curvatures of left invariant metrics on Lie groups, Adv.
in Math., 21, 293-329, 1976

[O1] Orta\c{c}gil, E.H., An Alternative Approach to Lie Groups and Geometric
Structures, Oxford University Press, 2018

[O2] Orta\c{c}gil, E.H., Curvature without connection, arXiv:2003.06593, 2020

[YB] Yano, K., Bochner, S., Curvature and Betti Numbers, Annals of Math.
Studies, No. 32, Princeton Univ. Press, 1953

\bigskip

Erc\"{u}ment Orta\c{c}gil

ortacgile@gmail.com

\end{document}